\documentclass[10pt,reqno,a4paper]{amsart}
\usepackage{amsfonts,amsmath,times}
\usepackage{a4wide}
\usepackage{tikz}

\usepackage[sc]{mathpazo} 
\usepackage[scaled]{helvet} 
\usepackage{eulervm} 

\definecolor{grau}{rgb}{0.65,0.65,0.65}
\definecolor{dblau}{rgb}{0,0,0.45}
\definecolor{blau}{rgb}{0,0,0.75} 
\definecolor{rot}{rgb}{0.9,0.1,0.1} 
\usepackage[pdftex]{hyperref}
\hypersetup{colorlinks,linkcolor=rot,citecolor=blue,urlcolor=blau}

\allowdisplaybreaks

\newtheorem{theorem}{Theorem}

\newtheorem{defi}{Definition}
\newtheorem{prop}{Proposition}
\newtheorem{coroll}{Corollary}
\theoremstyle{definition}
\newtheorem{remark}{Remark}
\newtheorem{example}{Example}

\newcommand{\myt}[1]{\emph{\color{dblau}#1}}

\newcommand{\JAP}{\emph{Journal of Applied Probability}}

\newcommand{\PTRF}{\emph{Probability Theory and Related Fields}}

\newcommand{\mom}{\text{model}\,\ensuremath{\mathcal{M}}}
\newcommand{\mor}{\text{model}\,\ensuremath{\mathcal{R}}}

\newcommand{\momp}{\text{model}\,\ensuremath{\mathcal{M}}\,}
\newcommand{\morp}{\text{model}\,\ensuremath{\mathcal{R}}\,}

\newcommand{\moms}{\text{model}\,\ensuremath{\mathcal{M}_{SEQ}}}
\newcommand{\mors}{\text{model}\,\ensuremath{\mathcal{R}_{SEQ}}}

\newcommand{\fallfak}[2]{\ensuremath{#1^{\underline{#2}}}}

\newcommand{\N}{\ensuremath{\mathbb{N}}}

\newcommand{\R}{\ensuremath{\mathbb{R}}}
\newcommand{\Z}{\ensuremath{\mathbb{Z}}}

\def\P{{\mathbb {P}}}
\def\E{{\mathbb {E}}}

\newcommand{\be}{\ensuremath{\boldsymbol{\mu}}}
\newcommand{\bs}{\ensuremath{\boldsymbol{\Sigma}}}

\newcommand{\folgen}{\ensuremath{\mathcal{S}}}

\newcommand\field{\mathcal{F}}

\newcommand\given{\, \vert \, }

\newcommand\Polya{P\'olya}

\author[Markus Kuba]{Markus Kuba}
\address{Markus Kuba\\
Institute of Applied Mathematics and Natural Sciences\\
University of Applied Sciences - Technikum Wien\\
H\"ochst\"adtplatz 5, 1200 Wien} 
\email{kuba@technikum-wien.at}

\title[Classification of urn models]{Classification of urn models with multiple drawings}
\keywords{Balanced urn model, multiple drawings, sampling without replacement, linear expected value}%
\subjclass[2000]{60C05, 60F05, 60G42}        
\begin{document}
\begin{abstract}
We consider multicolor urn models with multiple drawings.
An urn model is called linear if the conditional expected value of the urn composition at time $n$ is a linear function of the composition at time $n-1$. For four different sampling schemes - ordered and unordered samples with or without replacement - we classify urns into linear and non-linear models. We also discuss representations of the expected value and the covariance for linear models.
\end{abstract}

\date{\today}
\maketitle
\section{Introduction}
Urn models are simple, useful mathematical tools for describing many
evolutionary processes in diverse fields of application such as
analysis of algorithms and data structures, statistics and genetics.
 
The dynamics of the standard P{\'{o}}lya-Eggenberger urn models in the case of two types of colors and sample size $m=1$ can be described as follows. At the beginning, the urn contains $W_0$ white and $B_0$ black balls. At every step, we choose a
ball at random from the urn, examine its color and put it back into
the urn and then add/remove balls according to its color by the
following rules: if the ball is white, then we put $\alpha$ white and $\beta$
black balls into the urn, while if the ball is black, then $\gamma$ white
balls and $\delta$ black balls are put into the urn. The values $a, b, c,
d \in \mathbb{Z}$ are fixed integer values and the urn model is
specified by the ball transition matrix $M = \bigl(\begin{smallmatrix} a &  b \\
c & d \end{smallmatrix}\bigr)$. Models with $r\ge2$ types of
colors can be described in an analogous way and are specified by an
$r \times r$ ball transition matrix. One usually considers so-called \myt{tenable} urn models where the process of sampling and replacing balls never stops. Quantities of of interest are the number of white balls $W_n$ after~$n$ draws, and the number of black balls~$B_n$ after~$n$ draws in the case of $r=2$ colors; in the general case $r\ge 2$ one is interested in the distribution of the random vector $\mathbf{X}_n=(X^{(1)}_n,\dots, X^{(r)}_n)$, where $X_n^{(i)}$ denotes the number of balls colored $i$ after $n$ draws 
and the initial composition of the urn model at time $n=0$ is given by the (non-random) vector $\mathbf{X}_0=(X^{(1)}_0,\dots, X^{(r)}_0)$.  
In the classic version of \Polya\ urns with $r\ge 2$ different colors a single ball is sampled at each unit of discrete time.

Due to their importance in applications, there is a huge literature
on the stochastic behavior of urn models. The earliest contributions are the classical work of Ehrenfest and Ehrenfest~\cite{Ehrenfest}
and the article of Polya and Eggenberger~\cite{Eggenberger}.  We also refer to the classic surveys 
of Johnson and Kotz~\cite{JohnsonKotz1977,Kotz1997}, the book of Mahmoud~\cite{Mah2008} and the references therein.
The recent works of Chauvin et al.~\cite{Chauvin1,Chauvin2}, Janson~\cite{Jan2004,Jan2006,Jan2016},
Neininger and Knape~\cite{NeiningerKnape}, Pouyanne~\cite{Pou2008}, Mailler~\cite{Mailler}, M\"uller and Neininger~\cite{Mueller}, M\"uller~\cite{Mueller2}, are all devoted to urn models where only a \textit{single ball} is sampled at each step. 

\smallskip

In this work we are concerned with generalizations of so-called P{\'{o}}lya-Eggen\-berger urn
models. We study an $r$-color \Polya\ urn model, where \myt{multiple balls} are drawn at each discrete time step. 
We assume that the $r\ge 2$ different colors are in a fixed order and we thus speak of balls \myt{coloured} $i$, with $1\le i\le r$.
Their colors are inspected, then the sample is reinserted in the urn.

Additions and deletions take place according to the drawn sample of fixed size $m$ and we refer to the positive integer $m$ as the \myt{sample size};
 the sample is either a \myt{multiset} of size $m$ or a \myt{sequence} of length $m$. 
Such urn models recently received attention in the literature. 
Mahmoud and Tsukiji~\cite{TsukijiMahmoud2001} used an urn model with multiple drawings to study the distribution of random circuits,
Mahmoud studied urn models with sample size two~\cite{Mah2012}; we also refer the reader to the survery~\cite{JohnsonKotzMahmoud2004}.
Chen and Wei~\cite{ChenWei} generalized the original two-color Polya urn model to multiple drawings, see also Chen and Kuba~\cite{ChenKu2013+}.
A generalization of the Friedman urn model was an discussed in~\cite{KuMaPan2013+}. We also refer the reader to the general works of Moler et al.~\cite{Moler}, Higueras et al.~\cite{Higueras}, Renlund~\cite{Renlund},
as well as the recent very general work of Lasmar, Mailler and Selmi~\cite{Lasmar}.

\smallskip 

Urn models with multiple drawings and sample size $m>1$ are usually more difficult to analyze compared to the ordinary urn models with sample size $m=1$. 
The standard techniques - moment methods, analytic combinatorics and generating functions, embedding into continuous stochastic processes, and the contraction method  - are not easily applicable. In particular, the expected values of $X_n^{(i)}$ usually depend on the higher moments, making explicit and also asymptotic computations more complicated when the sample size is larger than one. 
It is however possible to obtain central limit theorems using the (Robbins–Monro) stochastic approximations techniques, see~\cite{Moler,Higueras,Mailler}, obtaining quite general limit theorems, also for unbalanced urn models.
A small drawback of these techniques is that expressions for the (positive integer) moments, in particular expectation and variance, moment convergence, as well as more precise information about the limit laws for so-called large-index urn models and also triangular urn models are at present elusive.
Thus, it is of interest to study classes of urn models with multiple drawings which generalize the existing very precise results for the standard case of sample size $m=1$.

\smallskip 

In the two-color case $\mathbf{X}_n=(X^{(1)}_n,X^{(2)}_n)=(W_n,B_n)$ with black and white balls it turned out~\cite{KuMa201314I} that a class of balanced urn models with multiple drawings and linear affine expected value contained the special cases treated before~\cite{ChenKu2013+,ChenWei,JohnsonKotzMahmoud2004,KuMaPan2013+,Mah2012,Moler,TsukijiMahmoud2001}.
This class was characterized by a condition on the first column entries of the 
ball replacement matrix: 
\begin{equation}
\label{MuliDrawsLinPropLinear}
a_k=(m-k)(a_{m-1}-a_m)+a_m,\qquad \mbox {for \ } 0\le k\le m,
\end{equation}
and the beforehand mentioned conditional expectation of the number of white balls $W_n$ after $n$ draws with an \myt{(affine) linear} structure of the form
\begin{equation}
\label{MuliDrawsLinPropLinearAlt}
\E\bigl[W_n \given \field_{n-1}\bigr]= \alpha_n W_{n-1} +\beta_n,\qquad n\ge 1.
\end{equation}
Here, $\alpha_n,\beta_n$ denote certain sequences depending only on the number of draws $n$, $a_{m-1}$, $a_m$, the total balance $\sigma$, and $\field_n$ 
denotes the sigma-algebra generated by the first $n$ draws from the urn:
\[
\alpha_n=\frac{T_{n-1}+m(a_{m-1}-a_m)}{T_{n-1}},\qquad \beta_n=a_m,\quad n\ge 1.
\]
For sample size $m=1$ the linear affine class reduces to ordinary balanced two-color urn models. The condition~\eqref{MuliDrawsLinPropLinear} and Equation~\eqref{MuliDrawsLinPropLinearAlt} allowed to obtain exact and asymptotic expressions for the expected value and the variance generalizing the previously known results for the case $m=1$, also leading to very precise 
limit laws~\cite{KuMa201314I,KuSu2016}, laws of the iterated logarithm. Moreover it was proven that thge martingale limits exhibit densities, bounded under suitable assumptions, and exponentially decaying tails.

\smallskip

We introduce a generalization of the condition~\eqref{MuliDrawsLinPropLinear} to the $r\ge 2$ color case for several different sampling models and the notion of \myt{linear} urn models with multiple drawings. 

\begin{defi}
\label{DefLinear}
A $r$-color urn model with multiple drawing and sample size $m\ge 1$ is called linear if the random vector $\mathbf{X}_n=(X^{(1)}_n,\dots, X^{(r)}_n)$, specifying the composition of the urn after $n$ draws, satisfies
\[
\E[\mathbf{X}_n \mid \field_{n-1}]=\mathbf{X}_{n-1}\cdot C_n, 
\]
for certain matrices $(C_k)_{k\in\N}\in\R^{r\times r}$.
\end{defi}

\begin{example}[Sample size $m=1$]
In the case of sample size $m=1$ all models are by definition linear. 
There, $C_k= I+\frac{1}{T_{k-1}}M$ and $M$ denoting the $r\times r$ ball transition matrix.
\end{example}

The definition above extends to unbalanced urn models.
There, the individual row sums of the ball replacement matrix $M$ are at least once different and
the total number of balls $T_n$ is itself a random variable. As a consequence, the matrices $(C_k)_{k\in\N}$ are random.
In this work we focus on the ball sampling point of view. However, it is possible to extend the definition 
to (generalized) Polya urn processes (compare with the discussions in~\cite{Jan2004,Jan2006,Pou2008}).

We consider four different sampling models and will determine all linear $r$-color urn models with multiple drawings and will show in Theorems~\ref{MDTheoremLinAff} and~\ref{MDTheoremLinAffSeq} that the resulting matrices $C_n$ are given by
\[
C_n=I+ \frac{1}{T_{n-1}}\cdot A,
\]
where $I$ denotes the $r\times r$ identity matrix and $A$ a certain $r\times r$ \myt{reduced ball transition matrix}, 
depending on the transition matrix $M$ of the urn model. 
We discuss the expected values $\be_n=\E[\mathbf{X}_n]$ and the (co)-variances $\boldsymbol\Sigma_n=\E[(\mathbf{X}_n-\bs)^T(\mathbf{X}_n-\bs)]$ of linear classes.
In the final section we briefly comment on limit laws. 

\subsection{Notation}
Throughout this work we use boldface letters $\mathbf{a}$, $\mathbf{k}$, etc. to denote row vectors. Given a vector $\mathbf{k}=k_1\mathbf{e}_1+\dots+k_r\mathbf{e}_r=(k_1,\dots,k_r)$ we frequently use the shorthand notation 
$\binom{m}{\mathbf{k}}=\binom{m}{k_1,\dots,k_r}$ for the multinomial coefficients. Moreover, given a vector $\mathbf{x}=(x_1,\dots,x_r)$ we use the notation $\mathbf{x}^\mathbf{k}:=x_1^{k_1}\cdot x_2^{k_2}\cdot\dots \cdot x_r^{k_r}$. We denote with $\fallfak{x}{k}$ the falling factorials $\fallfak{x}{k}=x(x-1)\dots(x-(k-1))$, and use the abbreviation
$\fallfak{\mathbf{x}}{\mathbf{k}}:=\fallfak{x_1}{k_1}\cdot\fallfak{x_2}{k_2} \dots \cdot \fallfak{x_r}{k_r}$.

\section{Preliminaries - Sampling schemes}
In the following we discuss in detail the four different sampling schemes and corresponding tenability assumptions. In the two-color case and $m=1$ tenability assumptions were classified in \cite{Bagchi1985}. In multiple drawings schemes, sufficient conditions for tenability under sampling without replacement, unordered samples, were formulated by Konzem and Mahmoud~\cite{KonzemMahmoud}.
\subsection{Sampling schemes: unordered samples}
\label{Sec:PrelimSampling}

We denote with $\{\big(X^{(1)}\big)^{k_1},\dots, \big(X^{(r)}\big)^{k_r}\}$ an unordered sample of size $m$ containing $k_i$ balls of color $i$, $1\le i\le r$. The vector $\mathbf{a}=(a_1,\dots,a_r)\in\Z^r$ with index $\mathbf{k}=k_1\mathbf{e}_1+\dots+k_r\mathbf{e}_r$ specifies that when drawing a multiset $\{\big(X^{(1)}\big)^{k_1},\dots, \big(X^{(r)}\big)^{k_r}\}$ of consisting of $k_i$ balls colored $i$, then we add/subtract $a_i$ balls colored $i$, $1\le i\le r$, with $\sum_{i=1}^r k_i=m$. The replacement of balls is specified by the the ball replacement matrix $M$ of dimension $\binom{r+m-1}{m-1}\times r$. It consists of $\binom{r+m-1}{m}$ row vectors $\mathbf{a}=\mathbf{a}_{\mathbf{k}}$. The \myt{index} $\mathbf{k}$ is contained in the discrete simplex 
\begin{equation}
\label{replacement2}
\Delta=\{\mathbf{k}=k_1\mathbf{e}_1+\dots+k_r\mathbf{e}_r =(k_1,\dots,k_r)\mid k_i\ge 0, \sum_{i=1}^{r}k_i=m\};
\end{equation}
here, $\mathbf{e}_i$ denotes the $i$th unit vector, such that
\begin{equation}
\label{replacement1}
M=(\mathbf{a}_\mathbf{k})_{\mathbf{k}\in\Delta}.
\end{equation}
Throughout this work we consider \myt{balanced} urn models (and briefly only comment on extensions to \myt{unbalanced} models) such that overall number of added/removed balls is a positive integer constant $\sigma>0$, independent of the composition of the sample;
Equivalently, we assume that the replacement matrix $M$ has constant row sum such that for all row vectors $\mathbf{a}_\mathbf{k}=(a_1,\dots,a_r)_\mathbf{k}$ of $M$ we have   
\begin{equation}
\label{balanced1}
\mathbf{a}_\mathbf{k} \cdot \boldsymbol{1}^T= \sum_{i=1}^{r}a_i=\sigma,
\end{equation}
with \myt{total balance} $\sigma>0$. Here $\boldsymbol{1}=\sum_{i=1}^{r}\mathbf{e}_i=(1,\dots,1)$.
As a consequence, the total number $T_n=\sum_{i=1}^{r}X_{n}^{(i)}$ of balls after $n$ draws is given by the \myt{deterministic} quantity
\begin{equation*}
T_n=n\cdot \sigma +T_0,\quad\text{with}\quad T_0=X^{(1)}_0+\dots+X^{(r)}_0=\mathbf{X}_0\cdot \boldsymbol{1}^T.
\end{equation*}

In the two-color case $r=2$ with sample size $r=2$ the ball replacement matrix $M$ is a $3\times 2$ matrix 
$M = 
\left(
\begin{smallmatrix}
\mathbf{a}_{2\cdot \mathbf{e}_1}\\
\mathbf{a}_{\mathbf{e}_1+ \mathbf{e}_2}\\
\mathbf{a}_{2\cdot \mathbf{e}_2}
\end{smallmatrix}
\right)
$; more commonly, the colors are addressed as white and black and $M$ is written as
$$
M = 
\bordermatrix{~ & {\color{grau}W }& {\color{grau}B }\cr
                  {\color{grau}\{WW\}} & a_0 & b_0 \cr
                  {\color{grau}\{WB\} }& a_1 & b_1 \cr
                 {\color{grau} \{BB\}} & a_2 & b_2 \cr
                  }.
$$

We consider two different sampling schemes for drawing an unordered sample of size $m$ at each step: \mom\ and \mor\footnote{The name of the models stems from the original works of Chen et al.~\cite{ChenKu2013+,ChenWei}}. In \mom\ we draw the $m$ balls without replacement. The $m$ balls are drawn at once and their colors are examined. After the sample is collected, we put the entire sample back in the urn and execute the replacement rules according to the counts of colors observed. The tenability assumption\footnote{These assumptions can be relaxed a little bit, if the initial values are adapted to the entries in the ball replacement matrix. E.g., for $m=1$ the urn model with ball replacement matrix
$\left(\begin{smallmatrix}-3 & 8\\
6& -4
\end{smallmatrix}\right)$ is still tenable because $W_0$ is a multiple of $3$ and $B_0$ a multiple of $4$.} for $0\le k\le m$. 
 implies that for \mom\ the coefficients $a_i$ of the vector
$\mathbf{a}=\mathbf{a}_{\mathbf{v}}=(a_1,\dots,a_r)_{\mathbf{v}}$, with index $\mathbf{k}=k_1\mathbf{e}_1+\dots+k_r\mathbf{e}_r$ satisfy 
\[
a_i\ge -k_i,\quad 1\le i\le r.
\]
We are never forced to remove more balls of color $i$ than previously drawn in the sample and the process of drawing and replacing balls can be continued ad infinitum.
Assume that an urn contains $c_i$ balls of color $i$, $1\le i\le r$, with $c_i\ge 0$. 
The probability $\P(\{\big(X^{(1)}\big)^{k_1},\dots, \big(X^{(r)}\big)^{k_r}\})$ of drawing $k_i$ balls of color $i$ is given by 
$$\P(\{\big(X^{(1)}\big)^{k_1},\dots, \big(X^{(r)}\big)^{k_r}\})= \binom{m}{k_1,\dots,k_r}\frac{\fallfak{c_1}{k_1}\dots\fallfak{c_r}{k_r}}{\fallfak{(c_1+\dots+c_r)}{m}}=\frac{\binom{c_1}{k_1}\dots\binom{c_r}{k_r}}{\binom{c_1+\dots+c_r}m},$$
with $k_i\ge 0$ and $\sum_{i=1}^r k_i=m$. The sample follows a multivariate hypergeometric distribution with $m$ draws and numbers $c_1,\dots,c_r$.

\smallskip

In \mor, we draw the $m$ balls with replacement. The $m$ balls are drawn one at a time. After a ball is drawn,  
its color is observed, and is reinserted in the urn, and thus 
it might reappear in the sampling of one multiset. After $m$ balls are collected in this way (and they are all back in the urn),
we execute the replacement rules according to the counts of colors observed.
By the tenability assumption the coefficients $a_i$ of the vector
$\mathbf{a}=\mathbf{a}_{\mathbf{k}}=(a_1,\dots,a_r)_{\mathbf{k}}$, with index $\mathbf{k}=k_1\mathbf{e}_1+\dots+k_r\mathbf{e}_r$ satisfy 
\[
a_i\ge -1,\quad \text{for} \quad k_i>0,\qquad \text{and} \quad a_i\ge 0 \quad\text{for} \quad k_i=0,
\]
$1\le i\le r$. The probability $\P(\{\big(X^{(1)}\big)^{k_1},\dots, \big(X^{(r)}\big)^{k_r}\})$ of drawing $k_i$ balls of color $i$ is given by 
$$\P(\{\big(X^{(1)}\big)^{k_1},\dots, \big(X^{(r)}\big)^{k_r}\})= \binom{m}{k_1,\dots,k_r}\frac{c_1^{k_1}\dots c_r^{k_r}}{(c_1+\dots+c_r)^{m}},$$
with $k_i\ge 0$ and $\sum_{i=1}^r k_i=m$.  Thus, the sample follows a multinomial distribution with $m$ trials and probabilities $p_i=c_i/(c_1+\dots+c_r)$, $1\le i\le m$.

Conditioning on the outcome of the $n$th draw, 
we obtain a distributional equation for the random vector $\mathbf{X}_n$.
The number of balls after $n$ draws is the number of  balls after $n-1$ draws, plus the contribution of $n$th draw:
\begin{equation}
\label{MDDistributional}
\mathbf{X}_n=\mathbf{X}_{n-1} + \sum_{\mathbf{k}\in \Delta} \mathbf{a}_{\mathbf{k}} \cdot \mathbb{I}_{n}\{\big(X^{(1)}\big)^{k_1},\dots, \big(X^{(r)}\big)^{k_r}\},
\end{equation}
for $n\ge 1$ with (non-random) initial composition $\mathbf{X}_0$. 
The indicators variables satisfy
\begin{equation*}
\P\bigl(\mathbb{I}_{n}\{\big(X^{(1)}\big)^{k_1},\dots, \big(X^{(r)}\big)^{k_r}\}=1\given \field_{n-1}\bigr)=\frac{\binom{X_{n-1}^{(1)}}{k_1}\dots \binom{X_{n-1}^{(r)}}{k_r} }{\binom{T_{n-1}}{m}}
\end{equation*}
for \mom, and 
\begin{equation*}
\P\bigl(\mathbb{I}_{n}\{\big(X^{(1)}\big)^{k_1},\dots, \big(X^{(r)}\big)^{k_r}\}=1\given \field_{n-1}\bigr)=\binom{m}{k_1,\dots,k_r}\frac{\big(X_{n-1}^{(1)}\big)^{k_1},\dots, \big(X_{n-1}^{(r)}\big)^{k_r}}{T_{n-1}^{m}}
\end{equation*}
for \mor.

\subsection{Sampling schemes: ordered samples}
\label{Sec:PrelimSamplingSeq}
We use the notation $\big(X^{(d_1)},\dots, X^{(d_m)}\big)$ to refer to an ordered sample of size $m$ such that at the $i$th draw for the ordered sample we got a ball colored $d_i$, $1\le d_i\le r$, $1\le i\le m$. The vector $\mathbf{a}=(a_1,\dots,a_r)\in\Z^r$ with index $\mathbf{d}=d_1\mathbf{e}_1+\dots+d_m\mathbf{e}_m$ specifies that when drawing $\big(X^{(d_1)},\dots, X^{(d_m)}\big)$ then we add/subtract $a_i$ balls colored $i$, $1\le i\le r$. We consider all possible sequences of outcomes of length $m$ :
$$\folgen=\{\mathbf{d}=(d_1,\dots,d_m) \mid 1\le d_i\le m,1\le i\le m\}.$$
The ball replacement matrix $M$ of the urn model is given by
\begin{equation}
\label{replacement3}
M=(\mathbf{a}_\mathbf{d})_{\mathbf{d}\in\folgen}
\end{equation}
and has dimension $r^m\times r$.
For example in the case $m=2$ and $r=2$ the ball replacement matrix $M$ is a $4\times 2$ matrix 
$\left(
\begin{smallmatrix}
\mathbf{a}_{(1,1)}\\
\mathbf{a}_{(1,2)}\\
\mathbf{a}_{(2,1)}\\
\mathbf{a}_{(2,2)}\\
\end{smallmatrix}
\right),
$
more commonly written as
$$
M = 
\bordermatrix{~ & {\color{grau} W} & {\color{grau}B}\cr
                  {\color{grau}(W,W)} & a_0 & b_0 \cr
                  {\color{grau}(W,B)} & a_1 & b_1 \cr
                  {\color{grau}(B,W) }& a_2 & b_2 \cr
                  {\color{grau}(B,B)} & a_3 & b_3 \cr
                  }.
$$

\smallskip

We consider again two different sampling schemes for drawing an ordered sample of size $m$ at each step: sampling without replacement in \moms\ and sampling with replacement in \mors. The tenability assumptions are identical to the cases of unordered samples. 
Given $\mathbf{d}\in\folgen$ let $\mathbf{j}(\mathbf{d})=(j_1(\mathbf{d}),\dots,j_r(\mathbf{d}))\in\Delta$, the discrete simplex $\Delta$ as defined in~\eqref{replacement2}, denote the vector of multiplicities:
\begin{equation}
\label{MuliDrawsLinDistEqnSeq0}
j_\ell(\mathbf{d})=\sum_{i=1}^{m}\mathbb{I}\{d_i=\ell\},\quad 1\le \ell \le r.
\end{equation}

Note that, given a replacement matrix $M=(\mathbf{a}_\mathbf{k})_{\mathbf{k}\in\Delta}$ associated to an unordered sampling scheme \momp or \morp,
there apparently exists an embedding into ordered sampling schemes. The corresponding replacement matrix $\hat{M}=(\hat{\mathbf{a}}_\mathbf{d})_{\mathbf{d}\in\folgen}$ can be defined as follows:
\[
\forall \mathbf{d}\in \folgen\colon \hat{\mathbf{a}}_{\mathbf{d}}=\mathbf{a}_{\mathbf{j}(\mathbf{d})}.
\]

Conditioning on the outcome of the $n$th draw, 
we obtain a distributional equation for the random vector $\mathbf{X}_n$. The indicators variables satisfy
\begin{equation*}
\P\bigl(\mathbb{I}_{n}\{\big(X^{(d_1)},\dots, X^{(d_m)}\big)\}=1\given \field_{n-1}\bigr)=\frac{\fallfak{(X_{n-1}^{(1)})}{j_1(\mathbf{d})}\dots \fallfak{(X_{n-1}^{(r)})}{j_r(\mathbf{d})} }{\fallfak{T_{n-1}}{m}}
\end{equation*}
for \moms, and 
\begin{equation*}
\P\bigl(\mathbb{I}_{n}\{\big(X^{(d_1)},\dots, X^{(d_m)}\big)\}=1\given \field_{n-1}\bigr)=\frac{\big(X_{n-1}^{(1)}\big)^{j_1(\mathbf{d})},\dots, \big(X_{n-1}^{(r)}\big)^{j_r(\mathbf{d})}}{T_{n-1}^{m}}
\end{equation*}
for \mors.
The distributional equation for the random vector $\mathbf{X}_n$ is given by
\begin{equation}
\label{MDDistributionalSeq}
\mathbf{X}_n=\mathbf{X}_{n-1} + \sum_{\mathbf{d}\in \folgen} \mathbf{a}_{\mathbf{d}} \cdot \mathbb{I}_{n}\{\big(X^{(d_1)},\dots, X^{(d_m)}\big)\}.
\end{equation}

\section{Classification of linear models}
\subsection{Unordered samples}
%

In the following we present a sufficient condition for (affine) \emph{linearity} of the conditional expected values of the $r$-color urn model with multiple drawings for both unordered samples and ordered samples.
We state first the result for unordered samples, generalizing the previous result of~\cite{KuMa201314I}.

\begin{theorem}
\label{MDTheoremLinAff}
A balanced $r$-color urn model with multiple drawings consisting of unordered samples of size $m$ is for both sampling schemes \momp and \morp linear if and only if the vectors $\mathbf{a}_{\mathbf{k}}$ of the ball replacement matrix $M=(\mathbf{a}_{\mathbf{k}})_{\mathbf{k}\in\Delta}$ are given by affine combinations of the vectors $\mathbf{a}_{m\cdot \mathbf{e}_i}$:
\[
\mathbf{a}_{\mathbf{k}}= \sum_{i=1}^{r}\frac{k_i}{m}\,\mathbf{a}_{m\cdot \mathbf{e}_i}.
\]
The conditional expected values is given by
\[
\E[\mathbf{X}_n \mid \field_{n-1}]=\mathbf{X}_{n-1}\cdot (I+ \frac{1}{T_{n-1}}\cdot A),
\]
with $A=(\mathbf{a}_{m\cdot \mathbf{e}_i})_{1\le i\le m}$ denoting the reduced $r\times r$ ball replacement matrix.
\end{theorem}

\begin{remark}
We assume above that $M$ is balanced with a certain total balance $\sigma>0$. 
Hence, we only can choose freely $r-1$ values in each row of the matrix with rows $\mathbf{a}_{m\cdot \mathbf{e}_i}$, $1\le i\le r$, 
due to the balance condition for the rows: $\mathbf{a}_{m\cdot \mathbf{e}_i} \cdot \boldsymbol{1}^T=\sigma$.
We call linear urn models \myt{triangular} if the reduced $r\times r$ ball replacement matrix $A$ is triangular.
\end{remark}

\begin{remark}
The condition above also naturally arises using an entirely different approach based on analytic combinatorics and generating functions. 
It turns out that the condition above ensures a system of ordinary differential equations for the expected value of $\mathbf{X}_n$
in contrast to the general case leading to (higher order) partial differential equations; this will be discussed elsewhere.
\end{remark}

\begin{example}[Two-color case]
In the two-color case $r=2$ and arbitrary $m\ge 1$ we re-obtain the previously derived condition for two colors white and black. Let $\mathbf{k}=(k,m-k)$. 
We get 
$$
\mathbf{a}_\mathbf{k}=\frac1m\big(k\mathbf{a}_{m\mathbf{e}_1}+(m-k)\mathbf{a}_{m\mathbf{e}_2}\big).
$$ 
Let $\mathbf{a}_{m\mathbf{e}_1}=(a_m,b_m)$, $\mathbf{a}_{m\mathbf{e}_2}=(a_0,b_0)$ and in general $\mathbf{a}_{\mathbf{k}}=(a_k,b_k)$, $0\le k\le m$.
Concerning the entries for the white balls we obtain the condition
\[
a_k=\frac1m \big(k a_m-(m-k)a_0\big)=a_0+k\cdot\frac{a_m-a_0}{m},\quad 0\le k\le m,
\] 
with $b_k=\sigma-a_k$. Rewriting the condition in terms of $a_{m-1}$ and $a_m$ instead of $a_0$ and $a_m$ we obtain
\[
a_k=(m-k)(a_{m-1}-a_m)+a_m, 
\]
as stated in~\cite{KuMa201314I}. The special case $m=2$ implies the condition of~\cite{Mah2012}.
\end{example}

\begin{example}[Generalized Polya urn model]
\label{ExamplePolya}
Chen and Wei~\cite{ChenWei} introduced a multicolor generalization of the classical Polya urn for \mom, which was then also considered under \mor~\cite{ChenKu2013+}: ``An urn contains balls of $r$ different
colors. For each color $i$, we initially have $X_0^{(i)}$ balls. At each step we draw $m$ balls at random, say $k_i$ balls
of color $i$, and their colors are noted. These balls are returned to the urn plus $c\cdot k_i$ balls of color $i$.''
According to this description we have for all $ \mathbf{k}\in\Delta\colon \mathbf{a}_{\mathbf{k}}=\sum_{i=1}^{m}c\cdot k_i \cdot \mathbf{e}_i$ and consequently the condition $\mathbf{a}_{\mathbf{k}}= \sum_{i=1}^{r}\frac{k_i}{m}\,\mathbf{a}_{m\cdot \mathbf{e}_i}$ is satisfied.
\end{example}

\begin{proof}
We obtain from the distributional equation the equation
\[
\E[\mathbf{X}_n \mid \mathcal{F}_{n-1}] = \mathbf{X}_{n-1} + \sum_{\mathbf{k}\in \Delta} \mathbf{a}_{\mathbf{k}}\cdot \E[\mathbb{I}_{n}\{\big(X^{(1)}\big)^{k_1},\dots, \big(X^{(r)}\big)^{k_r}\}].
\]
For \morp we get 
\begin{equation*}
\begin{split}
\E[\mathbf{X}_n \mid \mathcal{F}_{n-1}] 
&= \mathbf{X}_{n-1} + \sum_{\mathbf{k}\in \Delta} \mathbf{a}_{\mathbf{k}} \binom{m}{\mathbf{k}}\frac{\mathbf{X}_{n-1}^{\mathbf{k}}}{T_{n-1}^{m}}.
\end{split}
\end{equation*}
On the other hand the linearity implies that $\E[\mathbf{X}_n \mid \mathcal{F}_{n-1}]=\mathbf{X}_{n-1}\cdot C_n$, such that
\[
\mathbf{X}_{n-1}\cdot (C_n-I)= \sum_{\mathbf{k}\in \Delta} \mathbf{a}_{\mathbf{k}} \binom{m}{\mathbf{k}}\frac{\mathbf{X}_{n-1}^{\mathbf{k}}}{T_{n-1}^{m}}.
\]
Multiplication with $T_{n-1}^m$ gives
\[
\mathbf{X}_{n-1}\cdot T_{n-1}^{m}\cdot (C_n-I)=\sum_{\mathbf{k}\in \Delta} \mathbf{a}_{\mathbf{k}} \binom{m}{\mathbf{k}}\mathbf{X}_{n-1}^{\mathbf{k}}.
\]
In order to match powers we use the simple but important fact that 
\begin{equation}
T_{n-1}=\sum_{i=1}^{r}X_{n-1}^{(i)}.
\label{expansionR}
\end{equation}
We observe that unless $C_n=I+\frac{1}{T_n}B_n$, with $B_n\in\R^{r\times r}$, the powers of $\mathbf{X}_{n-1}^{\mathbf{k}}$ do not match on both sides of the equation.
Thus, 
\[
\mathbf{X}_{n-1}\cdot\Big(\sum_{i=1}^{r}X_{n-1}^{(i)}\Big)^{m-1}\cdot B_n=\sum_{\mathbf{k}\in \Delta} \mathbf{a}_{\mathbf{k}} \binom{m}{\mathbf{k}}\mathbf{X}_{n-1}^{\mathbf{k}}.
\]
Set $B_n=(\mathbf{b}_{n,1},\dots,\mathbf{b}_{n,r})$, we get further  
\[
\Big(\sum_{i=1}^{r}X_{n-1}^{(i)}\Big)^{m-1}\cdot \sum_{i=1}^{r}X_{n-1}^{(i)}\cdot \mathbf{b}_{n,i}  =\sum_{\mathbf{k}\in \Delta} \mathbf{a}_{\mathbf{k}} \binom{m}{\mathbf{k}}\mathbf{X}_{n-1}^{\mathbf{k}}.
\]
The multinomial theorem gives 
\[
\sum_{i=1}^{r} \sum_{\substack{\ell_1+\dots+\ell_{r}=m-1\\ \ell_j\ge 0}}\binom{m-1}{\mathbf{l}}
\mathbf{X}_{n-1}^{\mathbf{l}+\mathbf{e}_i}\cdot \mathbf{b}_{n,i}=\sum_{\mathbf{k}\in \Delta} \mathbf{a}_{\mathbf{k}} \binom{m}{\mathbf{k}}\mathbf{X}_{n-1}^{\mathbf{k}}.
\]
Comparison of coefficients of $\mathbf{X}_{n-1}^{\mathbf{k}}$ for every $\mathbf{k}\in \Delta$ gives the equation
\[
\sum_{i=1}^{r}\binom{m-1}{\mathbf{k}-\mathbf{e_i}} \mathbf{b}_{n,i} =\mathbf{a}_{\mathbf{k}} \binom{m}{\mathbf{k}},
\]
such that for $\mathbf{k}=m\cdot\mathbf{e}_i$ we obtain $\mathbf{b}_{n,i}= \mathbf{a}_{m\mathbf{e}_i}$ independent of $n$. 
Consequently, 
\[
\sum_{i=1}^{r}\binom{m-1}{\mathbf{k}-\mathbf{e_i}} \mathbf{a}_{m\mathbf{e}_i} =\mathbf{a}_{\mathbf{k}} \binom{m}{\mathbf{k}},
\]
and further
\[
\mathbf{a}_{\mathbf{k}}= \frac{1}{\binom{m}{\mathbf{k}}}\sum_{i=1}^{r}\binom{m-1}{\mathbf{k}-\mathbf{e_i}} \mathbf{a}_{m\cdot \mathbf{e}_i}.
\]
Since
\[
\frac{\binom{m-1}{\mathbf{k}-\mathbf{e_i}}}{\binom{m}{\mathbf{k}}}
=\frac{(m-1)!k_1!\dots k_i!\dots k_r!}{m!k_1!\dots k_{i-1}!(k_i-1)!k_{i+1}!\dots k_r!}=
\frac{k_i}{m},
\]
We obtain the stated result.

\smallskip 

For \momp we can can proceed in similar way. We get first
\begin{equation*}
\begin{split}
\E[\mathbf{X}_n \mid \mathcal{F}_{n-1}] &= \mathbf{X}_{n-1} + \sum_{\mathbf{k}\in \Delta} \mathbf{a}_{\mathbf{k}} \binom{m}{\mathbf{k}}
\frac{\fallfak{\mathbf{X}_{n-1}}{\mathbf{k}}}{\fallfak{T_{n-1}}{m}}.
\end{split}
\end{equation*}
Consequently, using the definition of the linearity
\[
\mathbf{X}_{n-1}\cdot\fallfak{T_{n-1}}{m}\cdot (C_n-I) = \sum_{\mathbf{k}\in \Delta} \mathbf{a}_{\mathbf{k}} \binom{m}{\mathbf{k}}
\fallfak{\mathbf{X}_{n-1}}{\mathbf{k}}.
\]
We use again the fact that $T_{n-1}=\sum_{i=1}^{r}X_{n-1}^{(i)}$ and the multinomial theorem 
for the falling factorials. We observe that unless 
$$
C_n-I=\Big(\frac{1}{T_{n-1}-j_1}\mathbf{b}_{n,1},\dots,\frac{1}{T_{n-1}-j_r}\mathbf{b}_{n,r}\Big),\quad 0 \le j_i\le m-1,
$$
the coefficients do not match since the left hand side has powers larger than on the right. 
Assume that for $1\le i\le r$ we have $j_i=0$ such that 
\begin{equation}
C_n-I=\frac{1}{T_{n-1}}\cdot B_n,
\label{AssumptionMOM}
\end{equation}
similar to \mor. Then, we can distribute the summand -1 in $\fallfak{(T_{n-1}-1)}{m-1}$ according to the value of $i$ to form $X_{n-1}^{(i)}-1$ and get
\begin{equation}
\label{expansionM}
\fallfak{(T_{n-1}-1)}{m-1}
= \sum_{\substack{\ell_1+\dots+\ell_{r}=m-1\\ \ell_j\ge 0}}\binom{m-1}{\mathbf{l}}\big((X_{n-1}^{(i)}-1)\big)^{\underline{\ell_i}}  \prod_{\substack{j=1\\j\neq i}}^{r}
\big(X_{n-1}^{(j)}\big)^{\underline{\ell_j}}
\end{equation}
valid for any $1\le i\le r$. Thus, we get further
\begin{equation}
\sum_{i=1}^{r} \sum_{\substack{\ell_1+\dots+\ell_{r}=m-1\\ \ell_j\ge 0}}\binom{m-1}{\mathbf{l}}
\fallfak{\mathbf{X}_{n-1}}{\mathbf{l}+\mathbf{e}_i}\cdot\mathbf{b}_{n,i}=\sum_{\mathbf{k}\in \Delta} \mathbf{a}_{\mathbf{k}} \binom{m}{\mathbf{k}}\fallfak{\mathbf{X}_{n-1}}{\mathbf{k}}.
\label{MDKoeffVergleich}
\end{equation}
Comparison of coefficients of $\fallfak{\mathbf{X}_{n-1}}{\mathbf{k}}$ for every $\mathbf{k}\in \Delta$ gives the equations 
\[
\sum_{i=1}^{r}\binom{m-1}{\mathbf{k}-\mathbf{e_i}} \mathbf{b}_{n,i} =\mathbf{a}_{\mathbf{k}} \binom{m}{\mathbf{k}},
\]
identical to before, leading to the stated result. 
It remains to justify~\eqref{AssumptionMOM}. Assume that $j_i$ is nonzero: we have
$$
\frac{\fallfak{T_{n-1}}{m}}{T_{n-1}-j_i}=\fallfak{T_{n-1}}{j_i}\cdot \fallfak{(T_{n-1}-j_i-1)}{m-1-j_i}.
$$
Using $\fallfak{T_{n-1}}{j_i}=\fallfak{(T_{n-1}-1)}{j_i}+j_i\fallfak{T_{n-1}}{j_i-1}$ we get
\[
\frac{\fallfak{T_{n-1}}{m}}{T_{n-1}-j_i}= \fallfak{(T_{n-1}-1)}{m}+j_i\cdot \fallfak{T_{n-1}}{j_i-1}\fallfak{(T_{n-1}-j_i-1)}{m-1-j_i}.
\]
From the first summands we can argue exactly the same way as for $j_i=0$ leading to $\mathbf{a}_{\mathbf{k}}= \sum_{i=1}^{r}\frac{k_i}{m}\,\mathbf{a}_{m\cdot \mathbf{e}_i}$. 
We observe that an additional contribution appears due to the second summand leading to additional powers in~\eqref{MDKoeffVergleich} non-existent on the right hand side. This is a contradiction, such that $j_i=0$.
Alternatively, by the properties of the multivariate hypergeometric distribution we get
\[
\E[\mathbf{X}_n \mid \mathcal{F}_{n-1}]=\mathbf{X}_{n-1} + \sum_{\mathbf{k}\in \Delta} \mathbf{a}_{\mathbf{k}} \binom{m}{\mathbf{k}}
\frac{\fallfak{\mathbf{X}_{n-1}}{\mathbf{k}}}{\fallfak{T_{n-1}}{m}}
=\mathbf{X}_{n-1} + \sum_{i=1}^{r}\frac{X_{n-1}^{(i)}}{T_{n-1}}\mathbf{a}_{m\cdot \mathbf{e}_i},
\]
such that we have a second rationale for $j_i=0$. 
\end{proof}

\subsection{Ordered samples}
For ordered samples we obtain the following counterpart of Theorem~\ref{MDTheoremLinAff}.

\begin{theorem}
\label{MDTheoremLinAffSeq}
A balanced $r$-color urn model with multiple drawings consisting of ordered samples of size $m$ is linear 
if and only if the vectors $\mathbf{a}_{\mathbf{d}}$ of $M$ satisfy for all $\mathbf{k}\in\Delta$ the condition:
\[
\sum_{\substack{\mathbf{d}\in\folgen\\\mathbf{j}(\mathbf{d})=\mathbf{k}}}\mathbf{a}_{\mathbf{d}}= \sum_{i=1}^{r}\binom{m-1}{\mathbf{k}-\mathbf{e_i}} \mathbf{a}_{i\cdot\boldsymbol{1}}.
\]
Then for both sampling schemes with or without replacement the identity 
\[
\E[\mathbf{X}_n \mid \mathcal{F}_{n-1}]=\mathbf{X}_{n-1}\cdot(I+ \frac{1}{T_{n-1}}\cdot A)
\]
holds, with $A=(\mathbf{a}_{i\cdot \boldsymbol{1}})_{1\le i\le m}$ denoting the reduced $r\times r$ matrix
and $\boldsymbol{1}=\sum_{j=1}^{m}\mathbf{e}_j$.
\end{theorem}
\begin{example}[Two-color case]
We consider the balanced two-color case $r=2$ with sample size $m=2$.
Let $\mathbf{a}_{1\cdot \mathbf{1}}=\mathbf{a}_{(1,1)}=(a_0,b_0)$, $\mathbf{a}_{(1,2)}=(a_1,b_1)$, $\mathbf{a}_{(2,1)}=(a_2,b_2)$ and $\mathbf{a}_{2\cdot \mathbf{1}}=\mathbf{a}_{(2,2)}=(a_3,b_3)$. 
We get the balancing condition $a_1+a_2=a_0+a_3$ or equivalently $b_1+b_2=b_0+b_3$.
\end{example}

\begin{example}[Embedding of $m=1$]
Given a quadratic replacement matrix $C\in\R^{r\times r}$ with row vectors $\mathbf{c}_i$, $1\le i \le r$, associated to sample size $m=1$. We can embed this case
using a suitable setup of the replacement matrix~\eqref{replacement3} $M=(\mathbf{a}_\mathbf{d})_{\mathbf{d}\in\folgen}$
such that $\forall\mathbf{d}\in\folgen:\ \mathbf{a}_{\mathbf{d}}\in\{\mathbf{c}_1,\dots,\mathbf{c}_r\}$. Given $\mathbf{k}\in\Delta$ we choose $\binom{m-1}{\mathbf{k}-\mathbf{e}_i}$ out of the $\binom{m}{\mathbf{k}}$ sequences $\mathbf{d}\in\folgen$ with $\mathbf{j}(\mathbf{d})=\mathbf{k}$ and define $\mathbf{a}_\mathbf{d}=\mathbf{c}_i$.
Consequently, $\mathbf{a}_{i\cdot\boldsymbol{1}}=\mathbf{c}_i$ and further
\[
\sum_{\substack{\mathbf{d}\in\folgen\\\mathbf{j}(\mathbf{d})=\mathbf{k}}}\mathbf{a}_{\mathbf{d}}=
\sum_{i=1}^{r}\binom{m-1}{\mathbf{k}-\mathbf{e_i}} \mathbf{c}_i
=\sum_{i=1}^{r}\binom{m-1}{\mathbf{k}-\mathbf{e_i}} \mathbf{a}_{i\cdot\boldsymbol{1}}.
\]
Due to the specific structure of $M$ the indicators can be grouped together and the distributional equation~\ref{MDDistributionalSeq} can be simplified mirroring 
essentially the case $m=1$.
\end{example}

\begin{proof}
Since the proof is similar to the proof of Theorem~\ref{MDTheoremLinAff} we will only discuss the main steps.
We obtain from~\eqref{MDDistributionalSeq} 
$$
\frac{1}{T_{n-1}}\mathbf{X}_n(C_n-I)=\sum_{\mathbf{d}\in \folgen} \mathbf{a}_{\mathbf{d}} \cdot \E[\mathbb{I}_{n}\{\big(X^{(d_1)},\dots, X^{(d_m)}\big)\}\mid\mathcal{F}_{n-1}].
$$
For \mors we observe that $C_n-I=\frac{1}{T_{n-1}}B_n$, with $B_n\in\R^{r\times r}$, such that
$$
T_{n-1}^{m-1}\mathbf{X}_nB_n=\sum_{\mathbf{d}\in \folgen} \mathbf{a}_{\mathbf{d}} \cdot \mathbf{X}_{n-1}^{\mathbf{j}(\mathbf{d})},
$$
with $\mathbf{j}(\mathbf{d})$ denoting the vector of multiplicities~\eqref{MuliDrawsLinDistEqnSeq0}. Expansion of $T_{n-1}$ as stated in~\eqref{expansionR} and extraction of coefficients from $\mathbf{X}_{n-1}^{\mathbf{k}}$, with $\mathbf{k}\in\Delta$, leads to the stated result. 
For \moms we observe again that $C_n-I=\frac{1}{T_{n-1}}B_n$, with $B_n\in\R^{r\times r}$, obtain an equation similar to the previous one:
$$
\fallfak{T_{n-1}}{m-1}\mathbf{X}_nB_n=\sum_{\mathbf{d}\in \mathcal{H}} \mathbf{a}_{\mathbf{d}} \cdot \fallfak{\mathbf{X}_{n-1}}{\mathbf{j}(\mathbf{d})}.
$$
We use~\eqref{expansionM} and compare coefficients from $\fallfak{\mathbf{X}_{n-1}}{\mathbf{k}}$, with $\mathbf{k}\in\Delta$, to get the stated result. 
\end{proof}

\section{Properties of linear models}
In the following we obtain the expected value and the covariance of linear models. First we turn to exact representations,
then we discuss asymptotic expansions.

\subsection{Expected value and covariance}
We collect a few properties of linear models.


We readily obtain a recurrence relation for the expected value from Theorem~\ref{MDTheoremLinAff}:
\[
\be_n=\be_{n-1}(I+\frac{1}{T_{n-1}}A),\quad n\ge 1,\qquad \be_0=\mathbf{X}_0.
\] 
Iteration of the recurrence relation gives the explicit result for $\be_n$.
\begin{prop}[Expected value for linear models]
For a linear $r$-color urn model with multiple drawings 
the expected value $\be_n=\E[\mathbf{X}_n]$ is given by
\[
\be_n=\mathbf{X}_0(I+\frac{1}{T_0}A)\dots(I+\frac{1}{T_{n-1}}A),
\]
where $A$ denotes the reduced ball replacement matrix of the linear model.
\end{prop}
Note that we can decompose the initial composition $\mathbf{X}_0$ in terms of a base of (generalized) left-eigenvectors of the ball replacement matrix $A$. Due to the balance condition the vector $(1,\dots ,1)$ is always an eigenvector corresponding to the largest eigenvalue $\sigma$. Specific assumptions on the algebraic multiplicity of the eigenvalues of $A$ readily lead to asymptotic expansions of the expected value.

\smallskip

The results for the covariance matrix $\bs_n=\E[(\mathbf{X}_n-\be_n)^T(\mathbf{X}_n-\be_n)]$ are not anymore model-independent in contrast to the expected value.
\begin{coroll}[Covariance matrix for linear models]
The covariance matrix $\bs_n=\E[(\mathbf{X}_n-\be_n)^T(\mathbf{X}_n-\be_n)]$ satisfies the following recurrence relations:
\begin{itemize}
	\item Unordered samples of size $m$:
	\begin{align*}
\bs_n&= (I+\frac{1}{T_{n-1}}A)^T\bs_{n-1}(I+\frac{1}{T_{n-1}}A)
-c_nA^T(\bs_{n-1}+\be_n^T\be_n)A\\
&\quad+\frac{1}{T_{n-1}}\sum_{i=1}^{r}\frac1m \E[X^{(i)}_{n-1}]\mathbf{a}_{m\cdot \mathbf{e}_i}^T\mathbf{a}_{m\cdot \mathbf{e}_i}
\end{align*}
with $c_n$ given by $c_n=\frac{1}{mT_{n-1}^2}$ for \momp and $c_n=\frac{1}{\fallfak{T_{n-1}}2}\Big(\frac1m-\frac{1}{T_{n-1}})$ for \mor.
\item Ordered samples of size $m$:
	\begin{align*}
\bs_n&= (I+\frac{1}{T_{n-1}}A)^T\bs_{n-1}(I+\frac{1}{T_{n-1}}A)
-\frac{1}{T_{n-1}^2}A^T\bs_{n-1}A+d_n\sum_{\mathbf{d}\in\folgen} \E[X^{\mathbf{j}(\mathbf{d})}]\mathbf{a}_{\mathbf{d}}^T\mathbf{a}_{\mathbf{d}},
\end{align*}
with $d_n$ given by $d_n=\frac{1}{\fallfak{T_{n-1}}{m}}$ for \moms\, and $d_n=\frac{1}{T_{n-1}^m}$ for \mors.
\end{itemize}
\end{coroll}
\begin{remark}
We observe that a linear model with unordered samples allows an exact and asymptotic computation of the covariance matrix. 
In contrast, for unordered samples the second (mixed) moments depend on higher moments $\E[X^{\mathbf{j}(\mathbf{d})}]$.
\end{remark}

\begin{proof}
Concerning the covariance matrix we study the mixed moments $\E[\mathbf{X}_n^T\mathbf{X}_n]$ and use then $\E[\mathbf{X}_n^T\mathbf{X}_n]=\bs-\be_n^T\be_n$.
The distributional equation~\eqref{MDDistributional} for $\mathbf{X}_n$ implies 
that 
\begin{align*}
\E[\mathbf{X}_n^T\mathbf{X}_n \mid \mathcal{F}_{n-1}] 
 &= \mathbf{X}_{n-1}^T\mathbf{X}_{n-1} + \sum_{\mathbf{k}\in\Delta}\Big(\mathbf{a}_{\mathbf{k}}^T\mathbf{X}_{n-1}+\mathbf{X}_{n-1}^T\mathbf{a}_{\mathbf{k}}+\mathbf{a}_{\mathbf{k}}^T\mathbf{a}_{\mathbf{k}}\Big)\\
 &\qquad\times \E[\mathbb{I}_{n}\{\big(X^{(1)}\big)^{k_1},\dots, \big(X^{(r)}\big)^{k_r}\}\mid \mathcal{F}_{n-1}] .
\end{align*}
For \morp we get
\begin{align*}
\E[\mathbf{X}_n^T\mathbf{X}_n \mid \mathcal{F}_{n-1}] 
 &= \mathbf{X}_{n-1}^T\mathbf{X}_{n-1} + \sum_{\mathbf{k}\in\Delta}\sum_{i=1}^r\frac{k_i}{m}\Big(\mathbf{a}_{m\mathbf{e}_i}^T\mathbf{X}_{n-1}+\mathbf{X}_{n-1}^T\mathbf{a}_{m\mathbf{e}_i}\Big)
 \binom{m}{\mathbf{k}}\frac{\mathbf{X}_{n-1}^{\mathbf{k}}}{T_{n-1}^{m}}\\
 &\quad  + \sum_{\mathbf{k}\in\Delta}\sum_{i=1}^r\sum_{j=1}^r\frac{k_ik_j}{m^2}\mathbf{a}_{m\mathbf{e}_i}^T\mathbf{a}_{m\mathbf{e}_j}\binom{m}{\mathbf{k}}\frac{\mathbf{X}_{n-1}^{\mathbf{k}}}{T_{n-1}^{m}}.
\end{align*}
By the properties of the multinomial distribution 
\[
 \sum_{\mathbf{k}\in\Delta}k_ik_j\binom{m}{\mathbf{k}}\frac{\mathbf{X}_{n-1}^{\mathbf{k}}}{T_{n-1}^{m}}
 =
 \begin{cases}
 m(m-1)\frac{X_{n-1}^{(i)}X_{n-1}^{(j)}}{T_{n-1}^2},\quad \text{for }i\neq j,\\
 m(m-1)\frac{\big(X_{n-1}^{(i)}\big)^2}{T_{n-1}^2}+m\frac{X_{n-1}^{(i)}}{T_{n-1}},\quad \text{for } i=j.
 \end{cases}
\]
Consequently, 
\begin{align*}
\E[\mathbf{X}_n^T\mathbf{X}_n \mid \mathcal{F}_{n-1}] 
 &= \mathbf{X}_{n-1}^T\mathbf{X}_{n-1} + \sum_{i=1}^r\frac{X_{n-1}^{(i)}}{T_{n-1}}\Big(\mathbf{a}_{m\mathbf{e}_i}^T \mathbf{X}_{n-1}+\mathbf{X}_{n-1}^T\mathbf{a}_{m\mathbf{e}_i}\Big)\\
  &\quad  + \frac{m-1}m \sum_{i=1}^r\sum_{j=1}^r\frac{X_{n-1}^{(i)}X_{n-1}^{(j)}}{T_{n-1}^2}\mathbf{a}_{m\mathbf{e}_i}^T\mathbf{a}_{m\mathbf{e}_j}
  +\frac1m \sum_{i=1}^r\frac{X_{n-1}^{(i)}}{T_{n-1}^2}\mathbf{a}_{m\mathbf{e}_i}^T\mathbf{a}_{m\mathbf{e}_j},
\end{align*}
and further
\begin{align*}
\E[\mathbf{X}_n^T\mathbf{X}_n \mid \mathcal{F}_{n-1}] 
 &= \mathbf{X}_{n-1}^T\mathbf{X}_{n-1} + \frac{1}{T_{n-1}}\Big(A^T\mathbf{X}_{n-1}^T\mathbf{X}_{n-1}+\mathbf{X}_{n-1}^T\mathbf{X}_{n-1}A\Big)\\
  &\quad  + \frac{m-1}{mT_{n-1}^2} A^T\mathbf{X}_{n-1}^T\mathbf{X}_{n-1}A
  +\frac1m \sum_{i=1}^r\frac{X_{n-1}^{(i)}}{T_{n-1}^2}\mathbf{a}_{m\mathbf{e}_i}^T\mathbf{a}_{m\mathbf{e}_j}.  
\end{align*}
We obtain the stated result for \morp with $C_n=\frac{1}{T_{n-1}^2}$ after rearranging the terms and using $\E[\mathbf{X}_n^T\mathbf{X}_n]=\bs-\be_n^T\be_n$ . 
For $\momp$ we proceed in an identical way. The only difference is that we use now the properties of the multivariate hypergeometric distribution
\[
 \sum_{\mathbf{k}\in\Delta}k_ik_j\binom{m}{\mathbf{k}}\frac{\fallfak{\mathbf{X}_{n-1}}{\mathbf{k}}}{\fallfak{T_{n-1}}{m}}
 =
 \begin{cases}
 m(m-1)\frac{X_{n-1}^{(i)}X_{n-1}^{(j)}}{\fallfak{T_{n-1}}2},\quad \text{for }i\neq j,\\
 m(m-1)\frac{\big(X_{n-1}^{(i)}\big)^2}{\fallfak{T_{n-1}}2}+m\frac{X_{n-1}^{(i)}}{T_{n-1}},\quad \text{for } i=j.
 \end{cases}
\]
For ordered samples we can proceed in a similar manner: from the distributional equation~\eqref{MDDistributionalSeq} for $\mathbf{X}_n$ we get
\begin{align*}
\E[\mathbf{X}_n^T\mathbf{X}_n \mid \mathcal{F}_{n-1}] 
 &= \mathbf{X}_{n-1}^T\mathbf{X}_{n-1} + \sum_{\mathbf{d}\in\folgen}\Big(\mathbf{a}_{\mathbf{d}}^T\mathbf{X}_{n-1}+\mathbf{X}_{n-1}^T\mathbf{a}_{\mathbf{d}}+\mathbf{a}_{\mathbf{d}}^T\mathbf{a}_{\mathbf{d}}\Big)\\
 &\qquad\times \E[\mathbb{I}_{n}\{\big(X^{(1)}\big)^{k_1},\dots, \big(X^{(r)}\big)^{k_r}\}\mid \mathcal{F}_{n-1}] .
\end{align*}
By linearity of the models we can simplify the sum $$\sum_{\mathbf{d}\in\folgen}\Big(\mathbf{a}_{\mathbf{d}}^T\mathbf{X}_{n-1}+\mathbf{X}_{n-1}^T\mathbf{a}_{\mathbf{d}}\Big)\E[\mathbb{I}_{n}\{\big(X^{(1)}\big)^{k_1},\dots, \big(X^{(r)}\big)^{k_r}\}\mid \mathcal{F}_{n-1}],$$ leading to the stated result.
\end{proof}

\section{Perspective and Acknowledgments}
In this note we obtained a classification of linear balanced multicolor urn models with multiple drawings and sample size $m$ greater or equal one. There are various directions for further investigations. It is highly likely that the results of Janson~\cite{Jan2016} for $m=1$ concerning convergence of the expected value and the covariances can be extended to the linear models. 
Concerning limit laws for model \morp the general results of Moler et al.~\cite{Moler} are applicable to prove a central limit theorem when the second largest eigenvalue of the reduced ball replacement matrix is less than $\frac{\sigma}2$; moreover, for the remaining models the more general results of~\cite{Lasmar} seem to ne applicable to all the other sampling models: \momp, \moms and \mors.
Concerning unordered samples and both \momp and \morp we expect that the algebraic approach of Pouyanne~\cite{Pou2008} for so-called \myt{large-index} and also \myt{triangular} urns can be suitably adapted to the linear urn models. Moreover, we believe that methods of~\cite{KuMa201314I,KuSu2016} can be adapted to the multicolor case, compare with the results of M\"uller~\cite{Mueller2} for sample size $m=1$. 

\smallskip

The author warmly thanks Hosam M.~Mahmoud for a great many valuable and interesting discussions concerning the generalization of the conditions in~\cite{Mah2012,KuMa201314I} leading to the present work.

\end{document}